\RequirePackage{ifpdf}
\ifpdf 
\documentclass[pdftex]{sigma}
\else
\documentclass{sigma}
\fi

\usepackage[matrix,arrow,curve]{xy}
\usepackage{rotating,array}

\begin{document}

\allowdisplaybreaks

\renewcommand{\thefootnote}{$\star$}

\renewcommand{\PaperNumber}{015}

\FirstPageHeading

\ShortArticleName{Universal Teichm\"uller Space}

\ArticleName{The Group of Quasisymmetric Homeomorphisms \\ of the
Circle and Quantization \\ of the Universal Teichm\"uller Space\footnote{This paper is a
contribution to the Special Issue on Kac--Moody Algebras and Applications. The
full collection is available at
\href{http://www.emis.de/journals/SIGMA/Kac-Moody_algebras.html}{http://www.emis.de/journals/SIGMA/Kac-Moody{\_}algebras.html}}}

\Author{Armen G. SERGEEV}

\AuthorNameForHeading{A.G. Sergeev}

\Address{Steklov Mathematical Institute, 8 Gubkina Str., 119991 Moscow, Russia}

\Email{\href{mailto:sergeev@mi.ras.ru}{sergeev@mi.ras.ru}}

\ArticleDates{Received July 29, 2008, in f\/inal form February 05,
2009; Published online February 08, 2009}

\Abstract{In the f\/irst part of the paper we describe the complex
geometry of the universal Teichm\"uller space $\mathcal T$, which
may be realized as an open subset in the complex Banach space of
holomorphic quadratic dif\/ferentials in the unit disc. The quotient
$\mathcal S$ of the dif\/feomorphism group of the circle modulo
M\"obius transformations may be treated as a smooth part of~$\mathcal T$. In the second part we consider the quantization of
universal Teichm\"uller space $\mathcal T$. We explain f\/irst how to
quantize the smooth part $\mathcal S$ by embedding it into a
Hilbert--Schmidt Siegel disc. This quantization method, however,
does not apply to the whole universal Teichm\"uller space~$\mathcal
T$, for its quantization we use an approach, due to Connes.}


\Classification{58E20; 53C28; 32L25}

\section{Introduction}

The universal Teichm\"uller space $\mathcal T$, introduced by
Ahlfors and Bers, plays a key role in the theory of quasiconformal
maps and Riemann surfaces. It can be def\/ined as the space of
quasisymmetric homeomorphisms of the unit circle $S^1$ (i.e.\
homeomorphisms of $S^1$, extending to quasiconformal maps of the
unit disc $\Delta$) modulo M\"obius transformations. The space
$\mathcal T$ has a natural complex structure, induced by its
realization as an open subset in the complex Banach space
$B_2(\Delta)$ of holomorphic quadratic dif\/ferentials in the unit
disc $\Delta$. The space $\mathcal T$ contains all classical
Teichm\"uller spaces $T(G)$, where $G$ is a Fuchsian group, as
complex submanifolds. The space $\mathcal
S:=\text{Dif\/f}_+(S^1)/\text{M\"ob}(S^1)$ of normalized
dif\/feomorphisms of the circle may be considered as a~``smooth" part
of $\mathcal T$.

Our motivation to study $\mathcal T$ comes from the string theory.
Physicists have noticed (cf.~\cite{Scherk,Bowick-Rajeev}) that the space
$\Omega_d:=C_0^\infty(S^1,\mathbb R^d)$ of smooth loops in the
$d$-dimensional vector space $\mathbb R^d$ may be identif\/ied with
the phase space of bosonic closed string theory. By looking at a
natural symplectic form $\omega$ on $\Omega_d$, induced by the
standard symplectic form (of type ``$dp\wedge dq$") on the phase
space, one sees that this form can be, in fact, extended to the
Sobolev completion of $\Omega_d$, coinciding with the space
$V_d:=H_0^{1/2}(S^1,\mathbb R^d)$ of half-dif\/ferentiable
vector-functions. Moreover, the latter space is the largest in the
scale of Sobolev spaces $H_0^s(S^1,\mathbb R^d)$, on which $\omega$
is correctly def\/ined. So the form $\omega$ itself chooses the
``right" space to be def\/ined on. From that point of view, it seems
more natural to consider $V_d$ as the phase space of bosonic string
theory, rather than $\Omega_d$. In this paper we set $d=1$ to
simplify the formulas and study the space
$V:=V_1=H_0^{1/2}(S^1,\mathbb R)$.

According to Nag--Sullivan \cite{Nag-Sullivan}, there is a
natural group, attached to the space $V{=}H_0^{1/2}\!(S^1,\!\mathbb R)$,
and this is precisely the group $\text{QS}(S^1)$ of quasisymmetric
homeomorphisms of the circle. Again one can say that the space $V$
itself chooses the ``right" group to be acted on. The group
$\text{QS}(S^1)$ acts on $V$ by reparametrization of loops and this
action is symplectic with respect to the form~$\omega$. The
universal Teichm\"uller space $\mathcal
T=\text{QS}(S^1)/\text{M\"ob}(S^1)$ can be identif\/ied by this action
with a~space of complex structures on $V$, compatible with $\omega$.

The second half of the paper is devoted to the quantization of the
universal Teichm\"uller space~$\mathcal T$. We start from the Dirac
quantization of the smooth part $\mathcal
S=\text{Dif\/f}_+(S^1)/\text{M\"ob}(S^1)$. This is achieved by
embedding of $\mathcal S$ into the Hilbert--Schmidt Siegel disc
$\mathcal D_{\text{HS}}$. Under this embedding the dif\/feomorphism
group $\text{Dif\/f}_+(S^1)$ is realized as a subgroup of the
Hilbert--Schmidt symplectic group $\text{Sp}_{\text{HS}}(V)$, acting
on the Siegel disc by operator fractional-linear transformations.
There is a holomorphic Fock bundle $\mathcal F$ over $\mathcal
D_{\text{HS}}$, provided with a projective action of
$\text{Sp}_{\text{HS}}(V)$, covering its action on $\mathcal
D_{\text{HS}}$. The inf\/initesimal version of this action is a
projective representation of the Hilbert--Schmidt symplectic Lie
algebra $\text{sp}_{\text{HS}}(V)$ in a f\/ibre $F_0$ of the Fock
bundle $\mathcal F$. This def\/ines the Dirac quantization of the
Siegel disc $\mathcal D_{\text{HS}}$. Its restriction to $\mathcal
S$ gives a projective representation of the Lie algebra
$\text{Vect}(S^1)$ of the group $\text{Dif\/f}_+(S^1)$ in the Fock
space $F_0$, which def\/ines the Dirac quantization of the space
$\mathcal S$.

However, the described quantization procedure does not apply to the
whole universal Teich\-m\"uller space $\mathcal T$. By this reason we
choose another approach to this problem, based on Connes
quantization. (We are grateful to Alain Connes for drawing our
attention to this approach, presented in \cite{Connes}.)
Brief\/ly, the idea is the following. The $\text{QS}(S^1)$-action on
$\mathcal T$, mentioned above, cannot be dif\/ferentiated in classical
sense (in particular, there is no Lie algebra, associated to
$\text{QS}(S^1)$). However, one can def\/ine a quantized inf\/initesimal
version of this action by associating with any quasisymmetric
homeomorphism $f\in\text{QS}(S^1)$ a quantum dif\/ferential $d^qf$,
being an integral operator on $V$ with kernel, given essentially by
the f\/inite-dif\/ference derivative of $f$. In these terms the
quantization of $\mathcal T$ is given by a representation of the
algebra of derivations of~$V$, generated by quantum dif\/ferentials
$d^qf$, in the Fock space $F_0$.



\pdfbookmark[1]{I. Universal Teichm\"uller space}{part1}
\section*{I.\ Universal Teichm\"uller space} \label{ch1}

\section[Group of quasisymmetric homeomorphisms of $S^1$]{Group of quasisymmetric homeomorphisms of $\boldsymbol{S^1}$}
\label{sec1}

\subsection[Definition of quasisymmetric homeomorphisms]{Def\/inition of quasisymmetric homeomorphisms}
\label{subsec11}

\begin{definition}
\label{def1} A homeomorphism $h:S^1\to S^1$ is called
\textit{quasisymmetric} if it can be extended to a quasiconformal
homeomorphism $w$ of the unit disc $\Delta$.
\end{definition}

Recall that a homeomorphism $w: \Delta\to w(\Delta)$, having locally
$L^1$-integrable derivatives (in generalized sense), is called
\textit{quasiconformal} if there exists a measurable complex-valued
function $\mu\in L^\infty(\Delta)$ with
$\Vert\mu\Vert_\infty:=\text{ess sup}_{z\in \Delta}|\mu(z)| =:k<1$
such that the following \textit{Beltrami equation}
\begin{equation}
\label{eq11} w_{\bar z}=\mu w_{z}
\end{equation}
holds for almost all $z\in\Delta$. The function $\mu$ is called a
\textit{Beltrami differential} or \textit{Beltrami potential} of $w$
and the constant $k$ is often indicated in the name of the
$k$-quasiconformal maps.

In the case when $k=0$ the homeomorphism $w$, satisfying
\eqref{eq11}, coincides with a conformal map from $D$ onto $w(D)$.
For a dif\/feomorphism $w$ its quasiconformality means that $w$
transforms inf\/initesimal circles into inf\/initesimal ellipses, whose
eccentricities (the ratio of the large axis to the small one) are
bounded by a common constant $K<\infty$, related to the above
constant $k=\Vert\mu\Vert_\infty$ by the formula
\[
K=\frac{1+k}{1-k}.
\]
The least possible constant $K$ is called the \textit{maximal
dilatation} of $w$ and is also sometimes indicated in the name of
$K$-quasiconformal maps.

The inverse of a quasiconformal map is again quasiconformal and the
same is true for the composition of quasiconformal maps. This
implies that orientation-preserving quasisymmetric homeomorphisms of
$S^1$ form a \textit{group of quasisymmetric homeomorphisms of the
circle} $\text{QS}(S^1)$ with respect to composition.

Any orientation-preserving dif\/feomorphism $h\in\text{Dif\/f}_+(S^1)$
extends to a dif\/feomorphism of the closed unit disc
$\overline{\Delta}$, which is evidently quasiconformal, according to
the above criterion. So $\text{Dif\/f}_+(S^1)\subset\text{QS}(S^1)$,
and we have the following chain of embeddings
\[
\text{M\"ob}(S^1)\subset\text{Dif\/f}_+(S^1)\subset\text{QS}(S^1)
\subset\text{Homeo}_+(S^1) .
\]
Here, $\text{M\"ob}(S^1)$ denotes the M\"obius group of
fractional-linear automorphisms of the unit disc $\Delta$,
restricted to $S^1$.

\subsection[Beurling-Ahlfors criterion]{Beurling--Ahlfors criterion}
\label{subsec12}

There is an intrinsic description of quasisymmetric homeomorphisms
of $S^1$ in terms of cross ratios. Recall that the \textit{cross
ratio} of four dif\/ferent points $z_1$, $z_2$, $z_3$, $z_4$ on the complex
plane is given by the quantity
\[
\rho=\rho(z_1,z_2,z_3,z_4):=
\frac{z_4-z_1}{z_4-z_2}:\frac{z_3-z_1}{z_3-z_2} .
\]
The equality of two cross ratios
$\rho(z_1,z_2,z_3,z_4)=\rho(\zeta_1,\zeta_2,\zeta_3,\zeta_4)$ is a
necessary and suf\/f\/icient condition for the existence of a
fractional-linear map of the complex plane, transforming the
quadruple $z_1$, $z_2$, $z_3$, $z_4$ into the quadruple
$\zeta_1$, $\zeta_2$, $\zeta_3$, $\zeta_4$. In the case of quasiconformal
maps the cross ratios of quadruples may change but in a controlled
way. This property, reformulated in the right way for
orientation-preserving homeomorphisms of $S^1$, yields a criterion
of quasisymmetricity, due to Ahlfors and Beurling.

The required property reads as follows: for an
orientation-preserving homeomorphism $h:S^1\to S^1$ it should exist
a constant $0<\epsilon<1$ such that the following inequality holds
\begin{equation}
\label{eq12} \frac12(1-\epsilon)\leq
\rho(h(z_1),h(z_2),h(z_3),h(z_4)) \leq \frac12(1+\epsilon)
\end{equation}
for any quadruple $z_1,z_2,z_3,z_4\in S^1$ with cross ratio
$\rho(z_1,z_2,z_3,z_4)=\frac12$.

\begin{theorem}[Beurling--Ahlfors, cf.~\cite{Ahlfors1966,Lehto}] \label{be-ahl} Suppose that $h:S^1\to S^1$ is
an orientation-preserving homeomorphism of $S^1$. Then it can be
extended to a quasiconformal homeomorphism $w:\Delta\to\Delta$ if
and only if it satisfies condition \eqref{eq12}.
\end{theorem}

Douady and Earle (cf.~\cite{Douady-Earle}) have found an
explicit extension operator $E$, assigning to a quasisymmetric
homeomorphism $h$ its extension to a quasiconformal homeomorphism
$w$ of $\Delta$, which is conformally invariant in the sense that
$g(w\circ h)=w\circ g(h)$ for any fractional-linear automorphism of
$\Delta$.

Though quasisymmetric homeomorphisms of $S^1$, in general, are not
smooth, they enjoy certain H\"older continuity, provided by the
following

\begin{theorem}[Mori, cf.~\cite{Ahlfors1966}]
\label{mori} Let $w:\Delta\to\Delta$ be a $K$-quasiconformal
homeomorphism of the unit disc onto itself, normalized by the
condition: $w(0)=0$. Then the following sharp estimate
\[
|w(z_1)-w(z_2)|<16|z_1-z_2|^{1/K}
\]
holds for any $z_1\neq z_2\in\Delta$. In other words, the
homeomorphism $w$ satisfies the H\"older condition of order $1/K$ in
the disc $\Delta$.
\end{theorem}

\section{Universal Teichm\"uller space}
\label{sec2}

\subsection[Definition of universal Teichm\"uller space]{Def\/inition of universal Teichm\"uller space}
\label{subsec21}

\begin{definition}
\label{def2} The quotient space
\[
\mathcal T:=\text{QS}(S^1)/\text{M\"ob}(S^1)
\]
is called the \textit{universal Teichm\"uller space}. It can be
identif\/ied with the space of \textit{normalized} quasisymmetric
homeomorphisms of $S^1$, f\/ixing the points $\pm1$ and $-i$.
\end{definition}

As we have pointed out earlier, there is an inclusion
\[
\text{Dif\/f}_+(S^1)/\text{M\"ob}(S^1)\hookrightarrow\mathcal
T=\text{QS}(S^1)/\text{M\"ob}(S^1) .
\]
We consider the homogeneous space
\[
\mathcal S:=\text{Dif\/f}_+(S^1)/\text{M\"ob}(S^1)
\]
as a ``smooth'' part of $\mathcal T$.

The space $\mathcal T$ can be provided with the
\textit{Teichm\"uller distance function}, def\/ined by
\begin{equation*}
\text{dist}(g,h)=\frac12\log K(h\circ g^{-1})
\end{equation*}
for any quasisymmetric homeomorphisms $g,h\in\mathcal T$, extended
to quasiconformal homeomorphisms of the disc $\Delta$. Here,
$K(h\circ g^{-1})$ denotes the maximal dilatation of the
quasiconformal map $h\circ g^{-1}$. This def\/inition does not depend
on the extensions of $g$, $h$ to $\Delta$ and def\/ines a metric on
$\mathcal T$. The universal Teichm\"uller space is a complete
connected contractible metric space with respect to the introduced
distance function (cf.~\cite{Lehto}). Unfortunately, this
metric is not compa\-tible with the group structure on $\mathcal T$,
given by composition of quasisymmetric homeomorphisms (cf.~\cite[Theorem~3.3]{Lehto}).

The term ``universal" in the name of the universal Teichm\"uller
space is due to the fact that $\mathcal T$ contains, as complex
submanifolds, all classical Teichm\"uller spaces $T(G)$, where $G$
is a Fuchsian group (cf.~\cite{Nag1988}). If a Riemann
surface $X$ is uniformized by the unit disc $\Delta$, so that
$X=\Delta/G$, then the corresponding Techm\"uller space $T(G)$ may
be identif\/ied with the quotient
\[
T(G)=\text{QS}(S^1)^G/\text{M\"ob}(S^1),
\]
where $\text{QS}(S^1)^G$ is the subset of $G$-invariant
quasisymmetric homeomorphisms in $\text{QS}(S^1)$. The universal
Teichm\"uller space $\mathcal T$ itself corresponds to the Fuchsian
group $G=\{1\}$.

Since quasisymmetric homeomorphisms of $S^1$ are def\/ined in terms of
quasiconformal maps of~$\Delta$, i.e.\ in terms of solutions of
Beltrami equation in $\Delta$, one can expect that there is a
def\/inition of $\mathcal T$ directly in terms of Beltrami
dif\/ferentials. Denote by $B(\Delta)$ the set of Beltrami
dif\/ferentials in the unit disc $\Delta$. It follows from above that
it can be identif\/ied (as a set) with the unit ball in the complex
Banach space~$L^\infty(\Delta)$.

Given a Beltrami dif\/ferential $\mu\in B(\Delta)$, we can extend it
to a Beltrami dif\/ferential $\check\mu$ on the extended complex plane
$\overline{\mathbb C}$ by setting $\check\mu$ equal to zero outside
the unit disc~$\Delta$. Then, applying the existence theorem for
quasiconformal maps on the extended complex plane $\overline{\mathbb
C}$ (cf.~\cite{Ahlfors1966}), we get a normalized
quasiconformal homeomorphism~$w^\mu$, satisfying Beltrami equation~\eqref{eq11} on~$\overline{\mathbb C}$ with potential~$\check\mu$.
This homeomorphism is conformal on the exterior $\Delta_-$ of the
closed unit disc~$\overline\Delta$ on~$\overline{\mathbb C}$ and
f\/ixes the points $\pm1$, $-i$. The image~$\Delta^\mu:=w^\mu(\Delta)$ of
$\Delta$ under the quasiconformal map $w^\mu$ is called a
\textit{quasidisc}. We associate with Beltrami dif\/ferential $\mu\in
B(\Delta)$ the normalized quasidisc~$\Delta^\mu$. Introduce an
equivalence relation between Beltrami dif\/ferentials in $\Delta$ by
saying that two Beltrami dif\/ferentials $\mu$ and $\nu$ are
equivalent if $w^\mu|_{\Delta_-}\equiv w^\nu|_{\Delta_-}$. Then the
universal Teichm\"uller space $\mathcal T$ will coincide with the
quotient
\[
\mathcal T = B(\Delta)/{\sim}
\]
of the space $B(\Delta)$ of Beltrami dif\/ferentials modulo introduced
equivalence relation. In other words, it coincides with \textit{the
space of normalized quasidiscs} in $\overline{\mathbb C}$.

\subsection{Complex structure of the universal Teichm\"uller space}
\label{subsec22}

We introduce a complex structure on the universal Teichm\"uller
space $\mathcal T$, using its embedding into the space of quadratic
dif\/ferentials.

Given an arbitrary point $[\mu]$ of $\mathcal T$, represented by a
normalized quasidisc $w^\mu(\Delta)$, consider a~map
\begin{equation*}
\mu\longmapsto S(w^\mu|_{\Delta_-}),
\end{equation*}
assigning to a Beltrami dif\/ferential $\mu\in[\mu]$ the Schwarz
derivative of the conformal map $w^\mu$ on~$\Delta$. Due to the
invariance of Schwarzian under M\"obius transformations, the image
of $\mu$ under the above map depends only on the class $[\mu]$ of
$\mu$ in $\mathcal T$. Moreover, it is a holomorphic quadratic
dif\/ferentials in $\Delta_-$. The latter fact follows from the
transformation properties of Beltrami dif\/ferentials, prescribed by
Beltrami equation (according to \eqref{eq11}, Beltrami dif\/ferential
behaves as a $(-1,1)$-dif\/ferential with respect to conformal changes
of variable). Composing the above map with a fractional-linear
biholomorphism of $\Delta_-$ onto the unit disc $\Delta$, we obtain
a map
\begin{equation*}
\Psi: \  \mathcal T\longrightarrow B_2(\Delta) ,\qquad
[\mu]\longmapsto \psi(\mu) ,
\end{equation*}
associating a holomorphic quadratic dif\/ferential $\psi(\mu)$ in
$\Delta$ with a point $[\mu]$ of the universal Teichm\"uller space
$\mathcal T$.

The space $B_2(\Delta)$ of holomorphic quadratic dif\/ferentials in
$\Delta$ is a complex Banach space, provided with a natural
hyperbolic norm, given by
\[
\Vert\psi\Vert_2:=\sup_{z\in\Delta}\big(1-|z|^2\big)^2|\psi(z)|
\]
for a quadratic dif\/ferential $\psi$. It can be proved (cf.~\cite{Lehto}) that $\Vert\psi[\mu]\Vert_2\leq 6$ for any
Beltrami dif\/ferential $\mu\in B(\Delta)$.

The constructed map $\Psi:\mathcal T\to B_2(\Delta)$, called a
\textit{Bers embedding}, is a homeomorphism of $\mathcal T$ onto an
open bounded connected contractible subset in $B_2(\Delta)$,
containing the ball of radius 1/2, centered at the origin (cf.~\cite{Lehto}).

Using the constructed embedding, we can introduce a complex
structure on the universal Teichm\"uller space $\mathcal T$ by
pulling it back from the complex Banach space $B_2(\Delta)$. It
provides $\mathcal T$ with the structure of a complex Banach
manifold. (Note that the topology on $\mathcal T$, induced by the
map $\Psi$, is equivalent to the one, determined by the
Teichm\"uller distance function.)

Moreover, the composition of the natural projection
\[
B(\Delta)\longrightarrow\mathcal T=B(\Delta)/{\sim}
\]
with the constructed map $\Psi$ yields a holomorphic map
\[
F: \ B(\Delta)\longrightarrow B_2(\Delta)
\]
with respect to the natural complex structure on $B(\Delta)$ (cf.~\cite{Nag1988}).

\pdfbookmark[1]{II. QS-action on the Sobolev space of half-differentiable
functions}{part2}

\section*{II.~QS-action on the Sobolev space of half-dif\/ferentiable
functions} \label{ch2}

\section[Sobolev space of half-differentiable functions on $S^1$]{Sobolev space of half-dif\/ferentiable functions on $\boldsymbol{S^1}$} \label{sec3}

\subsection[Definition]{Def\/inition}
\label{subsec31}

The \textit{Sobolev space of half-differentiable functions} on $S^1$
is a Hilbert space $V:=H_0^{1/2}(S^1,\mathbb R)$, consisting of
functions $f\in L^2(S^1,\mathbb R)$ with zero average over the
circle, having generalized derivatives of order $1/2$ again in
$L^2(S^1,\mathbb R)$. In terms of Fourier series, a function $f\in
L^2(S^1,\mathbb R)$ with Fourier series
\[
f(z)=\sum_{k\neq 0}f_kz^k ,\qquad f_k=\bar f_{-k} ,\qquad
z=e^{i\theta} ,
\]
belongs to $H_0^{1/2}(S^1,\mathbb R)$ if and only if it has a f\/inite
Sobolev norm of order $1/2$:
\begin{equation}
\label{eq31} \Vert f\Vert_{1/2}^2=\sum_{k\neq
0}|k||f_k|^2=2\sum_{k>0}k|f_k|^2<\infty .
\end{equation}

The space $H_0^{1/2}(S^1,\mathbb R)$ is well known and widely used
in classical function theory (cf.~\cite{Zygmund}). However,
our motivation to employ this space comes from its relation to
string theory (cf.\ below).

\subsection{K\"ahler structure}
\label{subsec32}

A symplectic form on $V$ is given by a 2-form $\omega:V\times
V\to\mathbb R$, def\/ined in terms of Fourier coef\/f\/icients of
$\xi,\eta\in V$ by
\begin{equation}
\label{eq32} \omega(\xi,\eta)=
2\,\text{Im}\,\sum_{k>0}k\xi_k\bar\eta_k .
\end{equation}
Because of \eqref{eq31}, this form is correctly def\/ined on $V$.
Moreover, $H_0^{1/2}(S^1,\mathbb R)$ is the largest Hilbert space in
the scale of Sobolev spaces $H_0^{s}(S^1,\mathbb R)$, $s\in\mathbb
R$, on which this form is def\/ined. It should be also underlined that
the form $\omega$ is the only natural symplectic form on $V$ (we
shall make this point clear in Section \ref{subsec41}).

We return to our motivation for studying the space $V$. It is well
known to physicists (cf., e.g.,
\cite{Scherk,Bowick-Rajeev}) that the space
$\Omega_d=C_0^\infty(S^1,\mathbb R^{d})$ of smooth loops in the
$d$-dimensional vector space~$\mathbb R^{d}$ can be identif\/ied with
the phase space of bosonic closed string theory. The space
$\Omega_d$ has a~natural symplectic form, which coincides with the
image of the standard symplectic form (of type ``$dp\wedge dq$") on
the phase space of closed string theory under the above
identif\/ication. This form, computed in terms of Fourier
decompositions, coincides precisely with the form $\omega$, given by~\eqref{eq32}. As we have remarked, the latter form may be extended
to the Sobolev space~$V_d:=H_0^{1/2}(S^1,\mathbb R^d)$ and this
space is the largest in the scale $H_0^s(S^1,\mathbb R^d)$ of
Sobolev spaces, on which $\omega$ is correctly def\/ined. One can say
that symplectic form $\omega$ ``chooses'' the Sobolev space~$V_d$.
This is in contrast to~$\Omega_d$, which was taken for the phase
space of string theory simply because it's easier to work with
smooth loops. By this reason, we f\/ind it more natural to consider~$V_d$ as the phase space of string theory, which motivates the study
of~$V_d$ in more detail. In our analysis we set $d=1$ for
simplicity.

Apart from symplectic form, the Sobolev space $V$ has a complex
structure $J^0$, which can be given in terms of Fourier
decompositions by the formula
\[
\xi(z)=\sum_{k\neq 0}\xi_kz^k\longmapsto
(J^0\xi)(z)=-i\sum_{k>0}\xi_kz^k+i\sum_{k<0}\xi_kz^k .
\]
This complex structure is compatible with symplectic form $\omega$
and, in particular, def\/ines a K\"ahler metric $g^0$ on $V$ by
$g^0(\xi,\eta):=\omega(\xi,J^0\eta)$ or, in terms of Fourier
decompositions,
\[
g^0(\xi,\eta)=2\text{Re\,}\sum_{k>0}k\xi_k\bar\eta_k .
\]
In other words, $V$ has the structure of a K\"ahler Hilbert space.

The complexif\/ication $V^\mathbb C=H_0^{1/2}(S^1,\mathbb C)$ of $V$
is a complex Hilbert space and the K\"ahler metric $g^0$ on $V$
extends to a Hermitian inner product on $V^\mathbb C$, given by
\begin{equation}
\label{eq33} \langle \xi,\eta\rangle =\sum_{k\neq 0}|k|\xi_k\bar\eta_k .
\end{equation}
We extend the symplectic form $\omega$ and complex structure
operator $J^0$ complex linearly to $V^\mathbb C$.

The space $V^\mathbb C$ is decomposed into the direct sum of the
form
\[
V^\mathbb C=W_+\oplus W_-  ,
\]
where $W_{\pm}$ is the $(\mp i)$-eigenspace of the operator
$J^0\in\text{End\,}V^\mathbb C$. In other words,
\[
W_+ =\Big\{f\in V^\mathbb C:\ f(z)=\sum_{k>0}f_kz^k\Big\}  ,\qquad
W_-=\overline W_+= \Big\{f\in V^\mathbb C:\ f(z)=\sum_{k<0}f_kz^k\Big\}  .
\]
The subspaces $W_\pm$ are isotropic with respect to symplectic form
$\omega$ and the splitting $V^\mathbb C=W_+\oplus W_-$ is an
orthogonal direct sum with respect to the Hermitian inner product
$\langle\cdot,\cdot\rangle$, given by~\eqref{eq33}.


\section[Grassmann realization of $\mathcal T$]{Grassmann realization of $\boldsymbol{\mathcal T}$} \label{sec4}


\subsection{QS-action on the Sobolev space} \label{subsec41}

Note that any homeomorphism $h$ of $S^1$, preserving the
orientation, acts on $L_0^2(S^1,\mathbb R)$ by change of variable.
In other words, there is an operator $T_h:L_0^2(S^1,\mathbb R)\to
L_0^2(S^1,\mathbb R)$, acting by
\[
T_h(\xi):=\xi\circ
h-\frac1{2\pi}\int_0^{2\pi}\xi\left(h(\theta)\right)d\theta .
\]
This operator has the following remarkable property.

\begin{proposition}[Nag--Sullivan \cite{Nag-Sullivan}]
\label{qs-act} The operator $T_h$ acts on $V$, i.e.\ $T_h:V\to V$, if
and only if $h\in\text{\rm QS}(S^1)$. Moreover, if $h$ extends to a
$K$-quasiconformal homeomorphism of the unit disc~$\Delta$, then the
operator norm of $T_h$ does not exceed $\sqrt{K+K^{-1}}$, where
$K=K(h)$ is the maximal dilatation of $h$.
\end{proposition}

Moreover, transformations $T_h$ with $h\in\text{QS}(S^1)$ generate
symplectic transformations of $V$.

\begin{proposition}[Nag--Sullivan \cite{Nag-Sullivan}]
\label{qs-symp} For any $h\in\text{\rm QS}(S^1)$ we have
\[
\omega(h^*(\xi),h^*(\eta))=\omega(\xi,\eta)
\]
for all $\xi,\eta\in V$. Moreover, the complex-linear extension of
$\text{\rm QS}$-action to the complexification~$V^{\mathbb C}$ preserves
the holomorphic subspace~$W_+$ if and only if
$h\in\text{\rm M\"ob}(S^1)$. In the latter case, $T_h$ acts as a unitary
operator on $W_+$.
\end{proposition}

We have pointed out in Section~\ref{subsec32} that the Sobolev
space $V$ is ``chosen" by the symplectic form $\omega$. In the same
way, one can say that the space $V$ chooses the reparametrization
group $QS(S^1)$. Indeed, this is the biggest reparametrization
group, leaving $V$ invariant, according to Proposition \ref{qs-act}.
On the other hand, it is a group of ``canonical transformations",
preserving the symplectic form $\omega$, according to Proposition
\ref{qs-symp}. So we have a natural phase space $(V,\omega)$
together with a natural group $QS(S^1)$ of its canonical
transformations.

Here is an assertion, making clear in what sense $\omega$ is a
unique natural symplectic form on $V$.

\begin{proposition}[Nag--Sullivan \cite{Nag-Sullivan}]
\label{nat-symp} Suppose that $\tilde\omega:V\times V\to\mathbb R$
is a continuous bilinear form on $V$ such that
\[
\tilde\omega(h^*(\xi),h^*(\eta))=\tilde\omega(\xi,\eta)
\]
for all $\xi,\eta\in V$ and all $h\in\text{\rm M\"ob}(S^1)$. Then
$\tilde\omega=\lambda\omega$ for some real constant $\lambda$. In
particular, $\tilde\omega$~is non-degenerate (if it is not
identically zero) and invariant under the whole group ${\rm QS}(S^1)$.
\end{proposition}

\subsection[Embedding of the universal Teichm\"uller space into
an infinite-dimensional Siegel disc]{Embedding of the universal Teichm\"uller space\\ into
an inf\/inite-dimensional Siegel disc} \label{subsec42}

The Propositions \ref{qs-act} and \ref{qs-symp} imply that
quasisymmetric homeomorphisms act on the Hilbert space~$V$ by
bounded symplectic operators. Hence, we have a map
\begin{equation}
\label{eq41} \mathcal
T=\text{QS}(S^1)/\text{M\"ob}(S^1)\longrightarrow
\text{Sp}(V)/\text{U}(W_+) .
\end{equation}
Here, $\text{Sp}(V)$ is the symplectic group of $V$, consisting of
linear bounded symplectic operators on $V$, and $\text{U}(W_+)$ is
its subgroup, consisting of unitary operators (i.e.\ the operators,
whose complex-linear extensions to $V^{\mathbb C}$ preserve the
subspace $W_+$).

In terms of the decomposition
\[
V^\mathbb C=W_+\oplus W_-
\]
any linear operator $A:V^\mathbb C\to V^\mathbb C$ is written in the
block form
\[
A=\begin{pmatrix} a & b\\ c & d
\end{pmatrix} .
\]
Such an operator belongs to symplectic group $\text{Sp}(V)$, if it
has the form
\[
A=\begin{pmatrix} a & b\\
\bar b & \bar a \end{pmatrix}
\]
with components, satisfying the relations
\begin{equation*}
\bar a^ta-b^t\bar b=1 ,\qquad \bar a^tb=b^t\bar a ,
\end{equation*}
where $a^t$, $b^t$ denote the transposed operators $a^t:W_-\to W_-$,
$b^t:W_-\to W_+$. The unitary group $\text{U}(W_+)$ is embedded into
$\text{Sp}(V)$ as a subgroup, consisting of diagonal block matrices
of the form
\[
A=\begin{pmatrix} a & 0\\ 0 & \bar a \end{pmatrix} .
\]

The space
\[
\text{Sp}(V)/\text{U}(W_+) ,
\]
standing on the right hand side of \eqref{eq41}, can be regarded as
an inf\/inite-dimensional analogue of the Siegel disc, since it may be
identif\/ied with the space of complex structures on~$V$, compatible
with $\omega$. Indeed, any such structure $J$ determines a
decomposition
\begin{equation}
\label{eq42*} V^\mathbb C =W\oplus\overline{W}
\end{equation}
of $V^\mathbb C$ into the direct sum of subspaces, isotropic with
respect to $\omega$. This decomposition is orthogonal with respect
to the K\"ahler metric $g_J$ on $V^\mathbb C$, determined by $J$ and
$\omega$. The subspaces~$W$ and~$\overline{W}$ are identif\/ied with
the $(-i)$- and $(+i)$-eigenspaces of the operator $J$ on $V^\mathbb
C$ respectively. Conversely, any decomposition \eqref{eq42*} of the
space $V^\mathbb C$ into the direct sum of isotropic subspaces
determines a complex structure $J$ on $V^\mathbb C$, which is equal
to $-iI$ on $W$ and $+iI$ on $\overline{W}$ and is compatible with
$\omega$. This argument shows that symplectic group $\text{Sp}(V)$
acts transitively on the space $\mathcal J(V)$ of complex structures~$J$ on $V$, compatible with~$\omega$. Moreover, a complex structure~$J$, obtained from a reference complex structure $J^0$ by the action
of an element~$A$ of~$\text{Sp}(V)$, is equivalent to $J^0$ if and
only if $A\in\text{U}(W_+)$. Hence,
\[
\text{Sp}(V)/\text{U}(W_+)=\mathcal J(V) .
\]

The space on the right can be, in its turn, identif\/ied with the
\textit{Siegel disc} $\mathcal D$, def\/ined as the set
\[
\mathcal D= \{Z:W_+\to W_-\text{\ is a symmetric bounded linear
operator with}\ \bar ZZ<I\} .
\]
The symmetricity of $Z$ means that $Z^t=Z$ and the condition $\bar
ZZ<I$ means that symmetric operator $I-\bar ZZ$ is positive
def\/inite. In order to identify $\mathcal J(V)$ with $\mathcal D$,
consider the action of the group $\text{Sp}(V)$ on $\mathcal D$,
given by fractional-linear transformations $A:\mathcal D\to\mathcal
D$ of the form
\[
Z\longmapsto (\bar aZ+\bar b)(bZ+a)^{-1} ,
\]
where $A=\begin{pmatrix} a & b\\ \bar b & \bar a
\end{pmatrix}\ \in\text{Sp}(V)$. The isotropy
subgroup at $Z=0$ coincides with the set of ope\-rators
$A\in\text{Sp}(V)$ such that $b=0$, i.e.\ with $\text{U}(W_+)$.

So the space
\[
\mathcal J(V)=\text{Sp}(V)/\text{U}(W_+)
\]
can be identif\/ied with the Siegel disc $\mathcal D$, and we have the
following

\begin{proposition}[Nag--Sullivan \cite{Nag-Sullivan}]
\label{uts-emb} The map
\[
\mathcal T=\text{QS}(S^1)/\text{\rm M\"ob}(S^1)\hookrightarrow \mathcal
J(V)=\text{\rm Sp}(V)/\text{U}(W_+)=\mathcal D
\]
is an equivariant holomorphic embedding of Banach manifolds.
\end{proposition}

For the smooth part $\mathcal S$ of the universal Teichm\"uller
space we can obtain a stronger version of this assertion by
replacing symplectic group $\text{Sp}(V)$ with its
\textit{Hilbert--Schmidt subgroup} $\text{Sp}_{\text{HS}}(V)$. By
def\/inition, this subgroup consists of bounded linear operators
$A\in\text{Sp}(V)$ with block representations
\[
A=\begin{pmatrix} a & b\\ \bar b & \bar a \end{pmatrix} ,
\]
in which the operator $b$ is Hilbert--Schmidt.

The map $f\mapsto T_f$ def\/ines an embedding
\[
\mathcal S\hookrightarrow \text{Sp}_{\text{HS}}(V)/\text{U}(W_+) .
\]
We identify, as above, the right hand side with a subspace $\mathcal
J_{\text{HS}}(V)$ of the space $\mathcal J(V)$ of compatible complex
structures on $V$. We call complex structures $J\in\mathcal
J_{\text{HS}}(V)$ \textit{Hilbert--Schmidt}. As before, the space
$\mathcal J_{\text{HS}}(V)$ of Hilbert--Schmidt complex structures
on $V$ can be realized as a~\textit{Hilbert--Schmidt Siegel disc}
\[
\mathcal D_{\text{HS}}=\{Z:W_+\to W_-\ \text{\,is a symmetric
Hilbert--Schmidt operator with $\bar ZZ<I$}\} .
\]
We have

\begin{proposition}[Nag \cite{Nag1992}]
\label{ds-emb} The map
\[
\mathcal S=\text{\rm Dif\/f}_+(S^1)/\text{\rm M\"ob}(S^1)\hookrightarrow
\mathcal
J_{{\rm HS}}(V)=\text{\rm Sp}_{{\rm HS}}(V)/\text{\rm U}(W_+)=\mathcal
D_{{\rm HS}}
\]
is an equivariant holomorphic embedding.
\end{proposition}

\pdfbookmark[1]{III.~Quantization of $\mathcal S$}{part3}

\section*{III.~Quantization of $\boldsymbol{\mathcal S}$} \label{ch3}

\section{Statement of the problem} \label{sec5}

\subsection{Dirac quantization}
\label{subsec51}

We start by recalling a general def\/inition of quantization of
f\/inite-dimensional classical systems, due to Dirac. A
\textit{classical system} is given by a pair $(M,\mathcal A)$, where
$M$ is the phase space and $\mathcal A$~is the algebra of
observables.

The \textit{phase space} $M$ is a smooth symplectic manifold of even
dimension $2n$, provided with a~symplectic 2-form $\omega$. Locally,
it is equivalent to the standard model, given by symplectic vector
space $M_0:=\mathbb R^{2n}$ together with standard symplectic form
$\omega_0$, given in canonical coordinates $(p_i,q_i)$,
$i=1,\dots,n$, on $\mathbb R^{2n}$ by
\[
\omega_0=\sum_{i=1}^n dp_i\wedge dq_i .
\]

The \textit{algebra of observables} $\mathcal A$ is a Lie subalgebra
of the Lie algebra $C^\infty(M,\mathbb R)$ of smooth real-valued
functions on the phase space $M$, provided with the Poisson bracket,
determined by symplectic 2-form $\omega$. In particular, in the case
of standard model $M_0=(\mathbb R^{2n},\omega_0)$ one can take for
$\mathcal A$ the \textit{Heisenberg algebra} $\text{heis}(\mathbb
R^{2n})$, which is the Lie algebra, generated by coordinate
functions $p_i,q_i$, $i=1,\dots,n$, and 1, satisfying the
commutation relations
\begin{gather*}
    \{p_i,p_j\} =\{q_i,q_j\}=0 ,\\
    \{p_i,q_j\} =\delta_{ij}\qquad\text{for} \quad i,j=1,\dots,n .
\end{gather*}

\begin{definition}
\label{def3} The \textit{Dirac quantization} of a classical system
$(M,\mathcal A)$ is an irreducible Lie-algebra representation
\[
r:\ \mathcal A\longrightarrow\text{End}^*H
\]
of the algebra of observables $\mathcal A$ in the algebra of linear
self-adjoint operators, acting on a complex Hilbert space $H$,
called the \textit{quantization space}. The algebra $\text{End}^*H$
is provided with the Lie bracket, given by the commutator of linear
operators of the form $\frac1{i}[A,B]$. In other words, it is
required that
\[
r\left(\{f,g\}\right)=\frac1{i}\left(r(f)r(g)-r(g)r(f)\right)
\]
for any $f,g\in\mathcal A$. We also assume the following
normalization condition: $r(1)=I$.
\end{definition}

For complexif\/ied algebras of observables $\mathcal A^\mathbb C$ or,
more generally, complex involutive Lie algeb\-ras of observables (i.e.\
Lie algebras with conjugation) their Dirac quantization is given by
an irreducible Lie-algebra representation
\[
r:\ \mathcal A^\mathbb C\longrightarrow\text{End\,}H\ ,
\]
satisfying the normalization condition and the conjugation law:
$r(\bar f)=r(f)^*$ for any $f\in\mathcal A$.

We are going to apply this def\/inition of quantization to
inf\/inite-dimensional classical systems, in which both the phase
space and algebra of observables are inf\/inite-dimensional. For
inf\/inite-dimensional algebras of observables it is more natural to
look for their projective Lie-algebra representations. The above
def\/inition of quantization will apply also to this case if one
replaces the original algebra of observables with its suitable
central extension.

\subsection{Statement of the problem}
\label{subsec52}

We start from the Dirac quantization of an inf\/inite-dimensional
system $(V,\mathcal A)$ with the phase space, given by the Sobolev
space of half-dif\/ferentiable functions $V:=H_0^{1/2}(S^1,\mathbb
R)$. The role of algebra of observables $\mathcal A$ will be played
by the semi-direct product
\[
\mathcal A=\text{heis}(V)\rtimes\text{sp}_{\text{HS}}(V) ,
\]
being the Lie algebra of the Lie group $\mathcal
G=\text{Heis}(V)\rtimes\text{Sp}_{\text{HS}}(V)$. The symplectic
Hilbert--Schmidt group $\text{Sp}_{\text{HS}}(V)$ was introduced in
Section~\ref{subsec32}, while the Heisenberg algebra $\text{heis}(V)$ and
the corresponding Heisenberg group $\text{Heis}(V)$ are def\/ined, as
in f\/inite-dimensional situation. Namely, the \textit{Heisenberg
algebra}~$\text{heis}(V)$ of $V$ is a central extension of the
Abelian Lie algebra~$V$, generated by coordinate functions. In other
words, it coincides, as a vector space, with
$\text{heis}(V)=V\oplus\mathbb R$, provided with the Lie bracket
\[
\left[(x,s),(y,t)\right]:=\left(0,\omega(x,y)\right) ,\qquad x,y\in
V,\quad s,t,\in\mathbb R .
\]
Respectively, the \textit{Heisenberg group} $\text{Heis}(V)$ is a
central extension of the Abelian group $V$, i.e.\ the direct product
$\text{Heis}(V)=V\times S^1$, provided with the group operation,
given by
\[
(x,\lambda)\cdot(y,\mu):=\big(x+y,\lambda\mu\,
e^{i\omega(x,y)}\big) .
\]

The choice of the introduced Lie algebra $\mathcal A$ for the
algebra of observables is motivated by the following physical
considerations. As we have pointed put, the space $V_d$ is a natural
Sobolev completion of the space $\Omega_d:=C_0^\infty(S^1,\mathbb
R^d)$ of smooth loops in $\mathbb R^d$. In the same way, the Lie
algebra $\mathcal A=\text{heis}(V)\rtimes\text{sp}_{\text{HS}}(V)$
is a natural extension of the Lie algebra
$\text{heis}(\Omega_d)\rtimes\text{Vect}(S^1)$, where
$\text{Vect}(S^1)$ is the Lie algebra of the dif\/feomorphism group
$\text{Dif\/f}_+(S^1)$. The algebra $\text{heis}(\Omega_d)$ can be
identif\/ied with the Lie algebra of coordinate functions on
$\Omega_d$, while the algebra $\text{Vect}(S^1)$ is generated by
certain quadratic functions on $\Omega_d$ (cf.~\cite{Bowick-Rajeev}). One can say that the Lie algebra
$\text{heis}(\Omega_d)\rtimes\text{Vect}(S^1)$ is an
inf\/inite-dimensional analogue of the Poincar\`e algebra of the
$d$-dimensional Minkowski space $M^d$, where $\text{heis}(\Omega_d)$
plays the role of the Lie algebra of translations of $M^d$, while
$\text{Vect}(S^1)$ is an analogue of the Lie algebra of hyperbolic
rotations of $M^d$.

\section{Heisenberg representation} \label{sec6}

In this Section we recall the well known Heisenberg representation
of the f\/irst component $\text{heis}(V)$ of algebra of observables
$\mathcal A$. A detailed exposition of this subject may be found in
\cite{Pressley-Segal,Kac-Raina,Berezin}.

\subsection{Fock space}
\label{subsec61}

Fix an admissible complex structure $J\in\mathcal J(V)$. It def\/ines
a polarization of $V$, i.e. a decomposition of $V^\mathbb C$ into
the direct sum
\begin{equation}
\label{eq61} V^\mathbb C=W\oplus\overline{W} ,
\end{equation}
where $W$ (resp.\ $\overline{W}$) is the $(-i)$-eigenspace (resp.\
$(+i)$-eigenspace) of the complex structure operator $J$. The
splitting \eqref{eq61} is the orthogonal direct sum with respect to
the Hermitian inner product $\langle z,w\rangle_J:=\omega(z,Jw)$, determined by
$J$ and sympletic form $\omega$.

The Fock space $F(V^\mathbb C,J)$ is the completion of the algebra
of symmetric polynomials on $W$ with respect to a natural norm,
generated by $\langle \cdot,\cdot\rangle_J$. In more detail, denote by $S(W)$ the
algebra of symmetric polynomials in variables $z\in W$ and introduce
an inner product on $S(W)$, def\/ined in the following way. It is
given on monomials of the same degree by the formula
\[
\langle z_1\cdot\dots\cdot z_n,z_1'\cdot\dots\cdot z_n'\rangle_J=
\sum_{\{i_1,\dots,i_n\}}\langle z_1,z'_{i_1}\rangle_J\cdot\dots\cdot\langle z_n,z'_{i_n}\rangle_J
,
\]
where the summation is taken over all permutations
$\{i_1,\dots,i_n\}$ of the set $\{1,\dots,n\}$ (the inner product of
monomials of dif\/ferent degrees is set to zero), and extended to the
whole algebra~$S(W)$ by linearity. The completion $\widehat{S(W)}$
of $S(W)$ with respect to the introduced norm is called the
\textit{Fock space} of $V^\mathbb C$ with respect to complex
structure $J$:
\[
F_J=F(V^\mathbb C,J):= \widehat{S(W)} .
\]

If $\{w_n\}$, $n=1,2,\dots$, is an orthonormal basis of $W$, then an
orthonormal basis of $F_J$ can be given by the family of polynomials
\begin{equation}
\label{fock-base}
P_K(z)=\frac{1}{\sqrt{k!}}\langle z,w_1\rangle_J^{k_1}\cdot\dots\cdot\langle z,w_n\rangle_J^{k_n}
,\qquad z\in W ,
\end{equation}
where $K=(k_1,\dots,k_n,0,\dots)$, $k_i\in\mathbb N\cup0$, and
$k!=k_1!\cdot\dots\cdot k_n!$.

\subsection{Heisenberg representation}
\label{subsec62}

There is an irreducible representation of the Heisenberg algebra
$\text{heis}(V)$ in the Fock space $F_J=F(V^\mathbb C,J)$, def\/ined
in the following way. Elements of $S(W)$ may be considered as
holomorphic functions on $\overline{W}$, if we identify $z\in W$
with a holomorphic function $\bar w\mapsto\langle w,z\rangle$ on
$\overline{W}$. Accordingly, $F_J$ may be considered as a subspace
of the space $\mathcal O(\overline{W})$ of functions, holomorphic on
$\overline{W}$. With this convention the \textit{Heisenberg
representation}
\[
r_J: \ \text{heis}(V)\longrightarrow\text{End}^*F_J
\]
of the Heisenberg algebra $\text{heis}(V)$ in the Fock space
$F_J=F(V^\mathbb C,J)$ is def\/ined by the formula
\begin{equation}
\label{heis-rep} r_J(v)f(\bar w):=-\partial_vf(\bar w)+\langle w,v\rangle_Jf(\bar
w) ,
\end{equation}
where $\partial_v$ is the derivative in direction of $v\in V$.
Extending $r_J$ to the complexif\/ied algebra $\text{heis}^\mathbb
C(V)$, we obtain
\[
r_J(\bar z)f(\bar w):=-\partial_{\bar z}f(\bar w)
\]
for $v=\bar z\in\overline{W}$ and
\[
r_J(z)f(\bar w):=\langle w,z\rangle_Jf(\bar w)
\]
for $z\in W$. We set also $r_J(c):=\lambda\cdot I$ for the central
element $c\in\text{heis}(V)$, where $\lambda$ is an arbitrary f\/ixed
non-zero constant.

Introduce the \textit{creation} and \textit{annihilation operators}
on $F_J$, def\/ined for $v\in V^\mathbb C$ by
\begin{equation*}
a_J^*(v):=\frac{r_J(v)-ir_J(Jv)}2 ,\qquad
a_J(v):=\frac{r_J(v)+ir_J(Jv)}2  .
\end{equation*}
In particular, for $z\in W$
\begin{equation*}
a_J^*(z)f(\bar w)=\langle w,z\rangle_Jf(\bar w)  ,
\end{equation*}
and for $\bar z\in\overline{W}$
\begin{equation*}
a_J(\bar z)f(\bar w)=-\partial_{\bar z}f(\bar w) .
\end{equation*}
For an orthonormal basis $\{w_n\}$ of $W$, we def\/ine the operators
\[
a^*_n:=a^*(w_n)  ,\qquad a_n:=a(\bar w_n)  ,\qquad n=1,2,\dots,
\]
and $a_0:=\lambda\cdot I$.

A vector $f_J\in F_J\setminus\{0\}$ is called the \textit{vacuum},
if $a_nf_J=0$ for $n=1,2,\dots$. In other words, it is a~non-zero
vector, annihilated by operators $a_n$. It is uniquely def\/ined by
$r_J$ (up to a~multiplicative constant) and in the case of the
initial Fock space $F_0=F(V,J^0)$ we set $f_0\equiv1$. Acting on
vacuum~$f_J$ by creation operators $a^*_n$, we can def\/ine the action
of representation $r_J$ on any polynomial, which implies the
irreducubility of $r_J$.

So we have the following

\begin{proposition}[cf.~\cite{Pressley-Segal,Kac-Raina,Berezin}]
\label{fock-rep} There is an irredicible Lie algebra representation
\[
r_J: \ \text{\rm heis}(V)\longrightarrow\text{\rm End}^*F_J
\]
of the Heisenberg algebra $\text{\rm heis}(V)$ in the Fock space
$F_J=F(V^\mathbb C,J)$, given by the formula~\eqref{heis-rep}.
\end{proposition}

We shall see in the next Section that this representation is
essentially unique.

\section{Symplectic group action on the Fock bundle} \label{sec7}

\subsection{Shale theorem}
\label{subsec71}

To construct an irreducible representation of the second component
$\text{sp}_{\text{HS}}(V)$ of the algebra of observables $\mathcal
A$, we study an action of the Hilbert--Schmidt symplectic group
$\text{Sp}_{\text{HS}}(V)$ on the Fock spaces $F_J$. This action is
provided by the following theorem.

\begin{theorem}[Shale]
\label{shale} The representations $r_0$ in $F_0$ and $r_J$ in $F_J$
are unitary equivalent if and only if $J\in\mathcal
J_{\text{\rm HS}}(V)$. In other words, for $J\in\mathcal
J_{\text{\rm HS}}(V)$ there exists a unitary intertwining operator $U_J:
F_0\to F_J$ such that
\[
r_J(v)=U_J\circ r_0(v)\circ U_J^{-1}.
\]
\end{theorem}

This theorem was proved by Shale \cite{Shale} in 1962, an
independent proof was given in Berezin's book
\cite{Berezin}, published in Russian in 1965 (Berezin
obtained also an explicit formula for the intertwining operator
$U_J$).

The following Proposition gives a description of $U_J$ in terms of
the Hilbert--Schmidt Siegel disc $\mathcal D_{\text{HS}}$, based on
the identif\/ication of $\mathcal J_{\text{HS}}(V)$ with $\mathcal
D_{\text{HS}}$.

\begin{proposition}[Segal \cite{Segal}]
\label{sp-act} There is a projective unitary action of the group
$\text{\rm Sp}_{\text{\rm HS}}(V)$ on Fock spaces, defined by the unitary
operator $U_J$, given by the formula \eqref{eq71} below.
\end{proposition}

Here is an idea of Segal's construction, details may be found in
\cite{Segal}. Given an admissible complex structure
$J\in\mathcal J_{\text{HS}}(V)$, we identify it with a point $Z$ in
the Siegel disc $\mathcal D_{\text{HS}}$. Regarding~$Z$ as an
element of the symmetric square $\widehat{ S}^2(W)$, we can
associate with it an element $e^{Z/2}$ of $\widehat{S(W)}$. The
inner product of two such elements has a simple expression
\[
\langle e^{Z_1/2},e^{Z_2/2}\rangle =\det(1-\bar Z_1Z_2)^{-1/2} .
\]
The normalized elements
\[
\epsilon_Z:=\det(1-\bar ZZ)^{1/4}e^{Z/2}
\]
play the role of \textit{coherent states} (cf., e.g.,
\cite{Berezin}). In terms of these states the action of the
group $\text{Sp}_{\text{HS}}(V)$ on Fock spaces, def\/ined by
\[
\text{Sp}_{\text{HS}}(V)\ni A=\begin{pmatrix} a & b\\
\bar b & \bar a \end{pmatrix}\longmapsto U_J: F_0\to
F_J\qquad\text{for}\quad J=A\cdot J^0 ,
\]
is given by the formula
\begin{equation}
\label{eq71} U_J: \ \epsilon_Z\longmapsto \mu\det(1+a^{-1}\bar
bZ)^{1/2}\epsilon_{A\cdot Z} ,
\end{equation}
where $\mu:\mathbb C^*\to S^1$ is the radial projection.

\subsection[Dirac quantization of $V$ and $\mathcal S$]{Dirac quantization of $\boldsymbol{V}$ and $\boldsymbol{\mathcal S}$}
\label{subsec72}

We can unite Fock spaces $F_J$ into a \textit{Fock bundle} over
$\mathcal D_{\text{HS}}$, having the following properties.

\begin{proposition}
\label{fock-bun} The Fock bundle
\[
\mathcal F:=\bigcup_{J\in\mathcal J(V)} F_J\longrightarrow\mathcal
J(V)=\mathcal D_{\text{\rm HS}}
\]
is a Hermitian holomorphic Hilbert space bundle over $\mathcal
D_{\text{\rm HS}}$. It can be provided with a projective unitary action
of the group $\text{\rm Sp}_{\text{\rm HS}}(V)$, covering the natural
$\text{\rm Sp}_{\text{\rm HS}}(V)$-action on the Siegel disc~$\mathcal
D_{\text{\rm HS}}$.
\end{proposition}

The proof of holomorphicity of the Fock bundle $\mathcal
F\to\mathcal D_{\text{HS}}$ is analogous to the proof of
holomorphicity of the determinant bundle over the Hilbert--Schmidt
Grassmannian, given in~\cite{Pressley-Segal}. Note that the
Fock bundle is trivial, since the Siegel disc $\mathcal
D_{\text{HS}}$ is contractible (even convex), so the statement
follows from the Hilbert space version of the Oka principle (cf.~\cite{Bungart}). An explicit trivialization of $\mathcal
F\to\mathcal D_{\text{HS}}$ is provided by the action~\eqref{eq71}.
This action def\/ines a projective unitary action of the group
$\text{Sp}_{\text{HS}}(V)$ on $\mathcal F$, covering the
$\text{Sp}_{\text{HS}}(V)$-action on Siegel disc $\mathcal
D_{\text{HS}}$.

The inf\/initesimal version of this action yields a projective
representation of the symplectic algebra $\text{sp}_{\text{HS}}(V)$
in the Fock space $F_0$. We present an explicit description of this
representation, due to Segal.

Recall that symplectic algebra $\text{sp}_{\text{HS}}(V)$ is the Lie
algebra of symplectic Hilbert--Schmidt group
$\text{Sp}_{\text{HS}}(V)$, which consists of linear operators $A$
in $V^\mathbb C$, having the following block representations
\[
A=\begin{pmatrix} \alpha & \beta\\ \bar\beta & \bar\alpha
\end{pmatrix} .
\]
Here, $\alpha$ is a bounded skew-Hermitian operator and $\beta$ is a
symmetric Hilbert--Schmidt operator on $F_0$. The complexif\/ied Lie
algebra $\text{sp}_{\text{HS}}(V)^{\mathbb C}$ consists of operators
of the form
\[
A=\begin{pmatrix} \alpha & \beta \\
         \bar\gamma & -\alpha^t
\end{pmatrix}  ,
\]
where $\alpha$ is a bounded operator, while $\beta$ and $\bar\gamma$
are symmetric Hilbert--Schmidt operators on $F_0$.

The projective representation of complexif\/ied symplectic algebra
$\text{sp}_{\text{HS}}(V)^{\mathbb C}$ is given by the formula
\begin{equation}
\label{sp-rep} \text{sp}_{\text{HS}}(V)^{\mathbb C}\ni
A=\begin{pmatrix} \alpha & \beta \\
                   \bar\gamma & -\alpha^t
          \end{pmatrix}
\longmapsto \rho(A)=D_\alpha+\frac12M_\beta+\frac12M^*_\gamma .
\end{equation}
Here, $D_\alpha$ is the derivation of $F_0$ in $\alpha$-direction,
def\/ined by
\[
D_\alpha f(\bar w)=\langle \alpha w,\partial_{\bar w}\rangle f(\bar w) .
\]
The operator $M_\beta$ is the multiplication operator on $F_0$,
def\/ined by
\[
M_\beta f(\bar w)=\langle \bar\beta w,\bar w\rangle f(\bar w) ,
\]
and the operator $M^*_\gamma$ is the adjoint of $M_\gamma$:
$M^*_\gamma f(\bar w)=\langle \gamma\partial_w,\partial_{\bar w}\rangle f(\bar
w)$.

This is a projective representation with cocycle
\begin{equation}
\label{eq72} [\rho(A_1),\rho(A_2)]-\rho([A_1,A_2])=
\frac12\,\text{tr}(\bar\gamma_2\beta_1-\bar\gamma_1\beta_2)I  ,
\end{equation}
intertwined with the Heisenberg representation $r_0$ of
$\text{heis}(V)$ in $F_0$.

Thus we have the following

\begin{proposition}[Segal \cite{Segal}]
\label{symp-rep} There is a projective unitary representation
\[
\rho:\ \text{\rm sp}_{\text{\rm HS}}(V)\longrightarrow\text{\rm End}^*F_0 ,
\]
given by formula \eqref{sp-rep} with cocycle \eqref{eq72}. This
representation intertwines with the Heisenberg representation $r_0$
of $\text{\rm heis}(V)$ in $F_0$.
\end{proposition}

The Heisenberg representation $r_0$ in the Fock space $F_0$,
described in Proposition \ref{fock-rep}, and symplectic
representation $\rho$, constructed in Proposition \ref{symp-rep},
def\/ine together Dirac quantization of the system
$(V,\widetilde{\mathcal A})$, where $\widetilde{\mathcal A}$ is the
central extension of $\mathcal A$, determined by \eqref{eq72}.

The constructed Lie-algebra representation of
$\text{sp}_{\text{HS}}(V)$ in the Fock space $F_0$ may be also
considered as Dirac quantization of a classical system, consisting
of the phase space $\mathcal
D_{\text{HS}}=\text{Sp}_{\text{HS}}(V)/\text{U}(W_+)$ and the
algebra of observables, given by the central extension of Lie
algebra $\text{sp}_{\text{HS}}(V)$.

The restriction of this construction to the smooth part $\mathcal
S=\text{Dif\/f}_+(S^1)/\text{M\"ob}(S^1)$ of the universal
Teichm\"uller space $\mathcal T=\text{QS}(S^1)/\text{M\"ob}(S^1)$
yields the Dirac quantization of $\mathcal S$. Namely, we have the
following

\begin{proposition}
\label{diff-fock-bun} The restriction of the Fock bundle $\mathcal
F\to\mathcal D_{\text{\rm HS}}$ to $\mathcal S$ is a Hermitian
holomorphic Hilbert space bundle
\[
\mathcal F:=\bigcup_{J\in\mathcal S} F_J\longrightarrow\mathcal S
\]
over $\mathcal S$. This bundle is provided with a projective unitary
action of the diffeomorphism group $\text{\rm Dif\/f}_+(S^1)$, covering
the natural $\text{\rm Dif\/f}_+(S^1)$-action on $\mathcal S$.
\end{proposition}

The $\text{Dif\/f}_+(S^1)$-action on the Fock bundle, mentioned in
Proposition, was explicitly constructed in
\cite{Goodman-Wallach}. The inf\/initesimal version of this
action yields a unitary projective representation of the Lie algebra
$\text{Vect}(S^1)$ in the Fock space $F_0$. We can consider this
construction as Dirac quantization of the phase space $\mathcal S$,
provided with the algebra of observables, given by the central
extension of the Lie algebra $\text{Vect}(S^1)$, called the
\textit{Virasoro algebra}.

\pdfbookmark[1]{IV. Quantization of $\mathcal T$}{part4}
\section*{IV.~Quantization of $\boldsymbol{\mathcal T}$} \label{ch4}

\section{Dirac versus Connes quantization} \label{sec8}

Unfortunately, the method, used in previous Chapter for the
quantization of $\mathcal S$, does not apply to the whole space
$\mathcal T$. Though we still can embed $\mathcal T$ into the Siegel
disc $\mathcal D$, we are not able to construct a projective action
of symplectic group $\text{Sp}(V)$ on Fock spaces. According to
theorem of Shale, it is possible only for the Hilbert--Schmidt
subgroup $\text{Sp}_{\text{HS}}(V)$ of $\text{Sp}(V)$. So one should
look for another way of quantizing the universal Teichm\"uller space
$\mathcal T$. We are going to use for that the ``quantized calculus''
of Connes and Sullivan, presented in Chapter~IV of the Connes' book~\cite{Connes} and \cite{Nag-Sullivan}.

Recall that in Dirac's approach we quantize a classical system
$(M,\mathcal A)$, consisting of the phase space $M$ and the algebra
of observables $\mathcal A$, which is a Lie algebra, consisting of
smooth functions on $M$. The quantization of this system is given by
a representation $r$ of $\mathcal A$ in a Hilbert space~$H$, sending
the Poisson bracket $\{f,g\}$ of functions $f,g\in\mathcal A$ into
the commutator $\frac1{i}[r(f),r(g)]$ of the corresponding
operators. In Connes' approach the algebra of observables $\mathfrak
A$ is an associative involutive algebra, provided with an exterior
dif\/ferential $d$. Its quantization is, by def\/inition, a~representation $\pi$ of $\mathfrak A$ in a Hilbert space $H$,
sending the dif\/ferential $df$ of a function $f\in\mathfrak A$ into
the commutator $[S,\pi(f)]$ of the operator $\pi(f)$ with a~self-adjoint symmetry operator $S$ with $S^2=I$. The dif\/ferential
here is understood in the sense of non-commutative geometry, i.e.\ as
a~linear map $d:\mathfrak A\to\Omega^1(\mathfrak A)$, satisfying the
Leibnitz rule (cf.~\cite{Connes}).

In the following table
we compare Connes and Dirac approaches to
quantization.

\begin{table}[h!]
\centering
\setbox1\hbox{$\mathcal A$ -- involutive Lie algebra}
\setbox2\hbox{$\mathfrak A$ -- involutive associative algebra}


\vspace{-1mm}

\setlength{\tabcolsep}{2mm}
\begin{tabular}{|>{\hfil}m{9mm}<{\hfil}|>{\hfil}m{\wd1}<{\hfil}|>{\hfil}m{\wd2}<{\hfil}|}%
\hline
\vphantom{$\Bigl|$}%
&Dirac approach&Connes approach
\\
\hline \vspace*{-19mm}{\begin{sideways}\parbox[c]{20mm}{Classical\\
system}\end{sideways}} & $(M,\mathcal A)$\ where: &
\vphantom{$\Big|$} $(M,\mathfrak A)$\ where:
\\
& $M$ -- phase space & $M$ -- phase space
\\
& $\mathcal A$ -- involutive Lie algebra & $\mathfrak A$ --
involutive associative
\\
& of observables & algebra of observables with
\\
& & dif\/ferential $d\colon\mathfrak A\to\Omega^1(\mathfrak A)$
\medskip
\\
\hline
\vspace*{-26.5mm}{\begin{sideways}\parbox[c]{22mm}{Quantization}\end{sideways}}
& Lie-algebra representation & representation
\vphantom{$\Bigl|$}%
\\
& $r\colon\mathcal A \to \text{End\,}H$, & $\pi\colon\mathfrak A
\to \text{End\,}H$,
\\
& sending & sending
\\
& $\{f,g\}\mapsto \frac1{i}[r(f),r(g)]$ & $df\mapsto[S,\pi(f)]$,
\\
& & where\ $S=S^*$, $S^2=I$
\medskip
\\
\hline
\end{tabular}\vspace{-1mm}
\end{table}

Reformulating the notion of Connes quantization of algebra of
observables $\mathfrak A$, one can say that it is a representation
of the algebra $\text{Der}(\mathfrak A)$ of derivations of
$\mathfrak A$ in the Lie algebra $\text{End\,}H$. Recall that a
\textit{derivation} of an algebra $\mathfrak A$ is a linear map:
$\mathfrak A\to\mathfrak A$, satisfying the Leibnitz rule. Clearly,
derivations of an algebra $\mathfrak A$ form a Lie algebra, since
the commutator of two derivations is again a derivation.

If all observables are smooth real-valued functions on $M$, the two
approaches are equivalent to each other. Indeed, the dif\/ferential
$df$ of a smooth function $f$ is symplectically dual to the
Hamiltonian vector f\/ield $X_f$ and this establishes a relation
between the associative algebra~$\mathfrak A$ of functions $f$ on
$M$ and the Lie algebra $\mathcal A$ of Hamiltonian vector f\/ields on
$M$. (This Lie algebra is isomorphic for a simply connected $M$ to a
Lie algebra of Hamiltonians, associated with $\mathcal A$.) A~symmetry operator $S$ is determined by a polarization $H=H_+\oplus
H_-$ of the quantization space~$H$. Evidently, $S=iJ$, where $J$ is
the complex structure operator, def\/ining the polarization
$H=H_+\oplus H_-$. (By this reason we do not make distinction
between symmetry and complex structure operators.)

In the case when the algebra of observables $\mathcal A$ contains
non-smooth functions, its Dirac quantization is not def\/ined in the
classical sense. In Connes approach the dif\/ferential $df$ of a~non-smooth observable $f\in\mathfrak A$ is also not def\/ined
classically, but its quantum counterpart $d^qf$, given by
\[
d^qf:=[S,\pi(f)]  ,
\]
may still be def\/ined, as it is demonstrated by the following
example, borrowed from~\cite{Connes}.

Suppose that $\mathfrak A$ is the algebra $L^\infty(S^1,\mathbb C)$
of bounded functions on the circle $S^1$. Any function
$f\in\mathfrak A$ def\/ines a bounded multiplication operator in the
Hilbert space $H=L^2(S^1,\mathbb C)$:
\[
M_f: \ v\in H\longmapsto fv\in H  .
\]
The operator $S$ is given by the \textit{Hilbert transform}
$S:L^2(S^1,\mathbb C)\to L^2(S^1,\mathbb C)$:
\[
(Sf)(e^{i\varphi})=\frac1{2\pi}\, V.P.\int_0^{2\pi}
K(\varphi,\psi)f(e^{i\psi})d\psi  ,
\]
where the integral is taken in the principal value sense and
$K(\varphi,\psi)$ is the \textit{Hilbert kernel}
\begin{equation}
\label{hil-ker} K(\varphi,\psi)= 1+i\cot\frac{\varphi-\psi}2 .
\end{equation}
The dif\/ferential $df$ of a general observable $f\in\mathfrak A$ is
not def\/ined in the classical sense, but its quantum analogue
\[
d^qf:=[S,M_f]
\]
is correctly def\/ined as an operator in $H$ for functions $f\in V$.
Namely, we have the following

\begin{proposition}[Nag--Sullivan~\cite{Nag-Sullivan}]
\label{hs-oper} A function $f\in V$ if and only if the corresponding
quantum differential $d^qf$ is a Hilbert--Schmidt operator on $H$
(and on $V$). Moreover, the Hilbert--Schmidt norm of $d^qf$
coincides with the $V$-norm of $f$.
\end{proposition}

Indeed, the commutator $d^qf:=[S,M_f]$ is an integral operator on
$H$ with the kernel, given by $K(\varphi,\psi)(f(\varphi)-f(\psi))$.
This operator is Hilbert--Schmidt if and only if its kernel is
square integrable on $S^1\times S^1$, i.e.{\samepage
\[
\int_0^{2\pi}\int_0^{2\pi}\frac{|f(\varphi)-f(\psi)|^2}
{\sin^2\frac12(\varphi-\psi)}\,d\varphi\,d\psi<\infty .
\]
This inequality is equivalent to the condition $f\in V$ (cf.~\cite{Nag-Sullivan}).}

The quantum dif\/ferential $d^qf=[S,M_f]$ of a function $f\in V$ is an
integral operator on $V$, given by
\begin{equation*}
d^qf(h)(e^{i\varphi})=\frac1{2\pi}\int_0^{2\pi}
k(\varphi,\psi)h(e^{i\psi})d\psi
\end{equation*}
with the kernel, given by
\[
k(\varphi,\psi):=K(\varphi,\psi)(f(\varphi)-f(\psi)) ,
\]
where $K(\varphi,\psi)$ is def\/ined by \eqref{hil-ker}.

Note that the quasiclassical limit of this operator, def\/ined by
taking the value of the kernel on the diagonal (i.e. by taking the
limit for $s\to t$), coincides (up to a constant) with the
multiplication operator $h\mapsto f'h$, so the quantization means in
this case essentially the replacement of the derivative by its
f\/inite-dif\/ference analogue. This f\/inite-dif\/ference analogue is an
integral operator, given by
\begin{equation}
\label{fin-diff} \delta f(v)(e^{i\varphi})=\frac1{2\pi}\int_0^{2\pi}
\frac{f(\varphi)-f(\psi)}{\varphi-\psi}v(e^{i\psi})d\psi .
\end{equation}

The correspondence between functions $f\in\mathfrak A$ and operators
$M_f$ on $H$ has the following remarkable properties (cf.~\cite{Power}):
\begin{enumerate}\itemsep=0pt
\item The dif\/ferential $d^qf$ is a f\/inite rank
operator if and only if $f$ is a rational function.
\item The dif\/ferential $d^qf$ is a compact operator if and only if
the function $f$ belongs to the class $\text{VMO}(S^1)$.
\item The dif\/ferential $d^qf$ is a bounded operator if and only if
the function $f$ belongs to the class $\text{BMO}(S^1)$.
\end{enumerate}
This list may be supplemented by further function-theoretic
properties of elements of $\mathfrak A$, having curious
operator-theoretic characterizations (cf.~\cite{Connes}).

\section{Quantization of the universal Teichm\"uller space}
\label{subsec91}

We apply these ideas to the universal Teichm\"uller space $\mathcal
T$. In Section \ref{subsec41} we have def\/ined a~natural action of
quasisymmetric homeomorphisms on $V$. As we have remarked, this
action does not admit the dif\/ferentiation, so classically there is
no Lie algebra, associated with $\text{QS}(S^1)$ or, in other words,
there is no classical algebra of observables, associated to
$\mathcal T$. (The situation is similar to that in the example
above.) We would like to def\/ine a quantum algebra of observables,
associated to $\mathcal T$.

First of all, extend the $\text{QS}(S^1)$-action on $V$ to symmetry
operators by setting
\begin{equation}
\label{eq92} S^h:=h\circ S\circ h^{-1}
\end{equation}
for $h\in\text{QS}(S^1)$. This action agrees with a natural action
of $\text{QS}(S^1)$ on the universal Teichm\"uller space $\mathcal
T=\text{QS}(S^1)/\text{M\"ob}(S^1)$, considered as a space of
compatible complex structures on $V$. The quantized inf\/initesimal
version of the action \eqref{eq92} is given by the integral operator
$d^qh:V\to V$, equal to $d^qh=[S,\delta h]$ with $\delta h$ given by
\eqref{fin-diff}.

Let us recall the steps of the Dirac quantization of Sobolev space
$V$:
\begin{enumerate}\itemsep=0pt
\item[1)] we start from $\text{Sp}_{\text{HS}}(V)$-action on $V$;
\item[2)] extend it to $\text{Sp}_{\text{HS}}(V)$-action on
complex structure operators $J$;
\item[3)] this action generates a projective unitary action of
$\text{Sp}_{\text{HS}}(V)$ on Fock spaces $F(V,J)$;
\item[4)] the inf\/initesimal version of this action yields a~projective unitary representation of the Lie algebra
$\text{sp}_{\text{HS}}(V)$ in Fock space $F_0$, described in
Section~\ref{subsec72}.
\end{enumerate}

In the case of $\mathcal T$ we have:
\begin{enumerate}\itemsep=0pt
\item[1)] $\text{QS}(S^1)$-action on $V$;
\item[2)] this action extends to $\text{QS}(S^1)$-action on symmetry
operators $S$, given by $h\mapsto S^h$.
\end{enumerate}
However, compared to Dirac quantization of $V$, the next step in the
quantization scheme is absent. Because of the Shale theorem, we
cannot extend the $\text{QS}(S^1)$-action on symmetry opera\-tors~$S$
to Fock spaces $F(V,S)$. Also we cannot dif\/ferentiate the
$\text{QS}(S^1)$-action on $V$. But we have a quantized
inf\/initesimal version of $h:S\mapsto S^h$, given by quantum
dif\/ferential $d^qh=[S,\delta h]$. We extend this operator to $F_0$
by def\/ining it f\/irst on the basis elements~\eqref{fock-base} of the
Fock space $F_0$ with the help of Leibnitz rule, and then by
linearity to all f\/inite elements of $F_0$. The completion of the
obtained operator yields an operator $d^qh$ on $F_0$. These extended
operators $d^qh$ with $h\in\text{QS}(S^1)$ generate a
\textit{quantum derivation algebra} $\text{Der}^q(\text{QS})$,
associated to $\mathcal T$. This algebra should be considered as a
quantum Lie algebra of observables, associated to $\mathcal T$. So,
instead of steps (3), (4) in the Dirac quantization of $V$, we
construct directly a quantum Lie algebra of observables
$\text{Der}^q(\text{QS})$, corresponding to the non-existing
classical Lie algebra of observables on~$\mathcal T$.

Moreover, we can use the quantum Lie algebra
$\text{Der}^q(\text{QS})$ as a substitution of a classical Lie
algebra of $\text{QS}(S^1)$.

\medskip

\noindent
\textbf{Conclusion.} The Connes quantization of the universal
Teichm\"uller space $\mathcal T$ consists of two stages:
\begin{enumerate}\itemsep=0pt
\item The f\/irst stage (``f\/irst quantization") is a construction of
quantized inf\/initesimal version of $\text{QS}(S^1)$-action on $V$,
given by quantum dif\/ferentials $d^qh=[S,\delta h]$ with
$h\in\text{QS}(S^1)$.
\item The second step (``second quantization") is an extension of
quantum dif\/ferentials $d^qh$ to the Fock space $F_0$. The extended
operators $d^qh$ with $h\in\text{QS}(S^1)$ generate the quantum
algebra of observables $\text{Der}^q(\text{QS})$, associated to
$\mathcal T$.
\end{enumerate}

We note also that the correspondence principle for the constructed
Connes quantization of~$\mathcal T$ means that this quantization
reduces to the Dirac quantization while restricted to $\mathcal S$.

\subsection*{Acknowledgements}
While preparing this paper, the author was partly supported by the
RFBR grants 06-02-04012, 08-01-00014, by the program of Support of
Scientif\/ic Schools (grant NSH-3224.2008.1), and by the Scientif\/ic
Program of RAS ``Nonlinear Dynamics''.

\pdfbookmark[1]{References}{ref}
\LastPageEnding

\end{document}